\medskip 
\medskip
{\bf THE CLASSIFICATION OF THREE-DIMENSIONAL GRADIENT LIKE MORSE-SMALE DYNAMIC 
SYSTEMS}  
\medskip
{\bf A. O. Prishlyak}
\bigskip
In papers [1 --- 4] the topological classification of Morse-Smale vector fields
on 2-manifolds and in [4,5] on 3-manifolds is got. In [6] classification 
three-dimensional Morse-Smale diffeomorphisms with non-intersected stable and 
unstable manifolds of saddle points is given.  

In this paper new approach to the classifications problem of three-dimensional 
gradient like Morse-Smale dynamic systems is represented. The criterion of the 
such systems topological equivalence is in terms of homeomorphism between 
surfaces with two circles series on its is given.  

The invariant of vector fields and diffeomorphisms are constructed and the 
classification of gradient like dynamic systems is obtained.  
\bigskip
{\bf 1. Basic definitions.} 
\medskip
Smooth dynamic system (vector field or 
diffeomorphism) is called Morse-Smale system if: 

\medskip
1) it has the finite number of the critical elements (periodical trajectories for 
diffeomorphism, fixed points and the closed orbits in case of vector field) and
all of them is non-degenerate (hyperbolic);  

2) stable and unstable integral manifolds of critical elements have transversal 
intersections;  

3) the limit set for every trajectory is critical element. 
\medskip
Heteroclinic trajectories of diffeomorphism are trajectories lying in the 
stable and unstable manifolds intersection for the same index critical 
elements. 

Morse-Smale dynamic system is called by gradient like, if there 
isn't closed orbit in case of vector field [7] and heteroclinic 
trajectories in the case of diffeomorphism [8]. 

Vector fields are called topological equivalent, if there exists the 
homeomorphism of manifold onto itself, which maps integral trajectories into 
integral trajectories preserving their orientation. By graph we will understand
finite 1-dimensional CW-complex. The isomorphism of graphs is the cell 
homeomorphism (id est., graph homeomorphism, which maps vertexes into vertexes 
and edges into edges).

\bigskip
{\bf 2. The criterion of the vector fields topological equivalence.}  
\medskip
Let $M^3$ be a closed oriented manifold, $X$ and $X'$ be Morse-Smale vector 
fields on it. Let $a_1,...$, $a_k$ be fixed 0-points of the fields $X$, and 
$a_1',...$, $a_k'$ of the fields $X'$; $b_1,...$, $b_n$ and $b_1',...$, $b_n'$ 
be the fixed 1-points. Let $K$ is the union of the 0- and 1-points stable 
manifolds. We shall consider tubular neighborhood $U (K)$ of this union.

We denote by $N= \partial U(K)$ the boundary of this neighborhood for field \~O
and by $N'$ for field $X'$ . Then these boundaries are the surfaces which 
Heegaard splitting of manifold $M^3 [9]$.  

We denote by $v(x)$ and $u(x)$ stable and unstable manifolds of fixed points 
$x$ . Let $u\__i$ be the circles, which are obtained as a result of 1-points 
unstable manifold intersections. Then $u\__i$ is a set of non crossed circles 
on surface $F$.  

If $c_1,..$, $c_m$ are the fixed 2-points, then intersection $v_i=v(c_i) \cap  
N$ will form another set of circles on surface F. Analogously, for the field 
$X'$ on the surface $N'$ there exists two sets of circles.  

If there is one 0-point and one 3-point then sets of circles will be the 
systems of the meridians of surface which form Heegaard diagram of manifold 
$M^3 [9].$ 
\medskip
{\bf Lemma 1.} Field $X$ is topological equivalently to $X'$ if and only if 
there is a homeomorphism of surface $f: F\rightarrow  F$, which maps the first 
set of circles into the first one and the second into the second. 
\medskip
Proof. Necessity results from building. We shall prove adequacy. Let such 
homeomorphism exists. We shall consider disks which lie onto unstable integral 
manifolds $U(b_i)$ , contain dots $b_i$ and bound by circles $u_i$ . Then we 
can continue homeomorphisms from the boundary of these disks up to disks 
homeomorphisms and such that translate integral trajectories into integral 
trajectories (because each integral trajectory, except fixed dots, crosses the 
boundary of disk.  

Analogously there exists disks homeomorphisms consisting of the of integral 
trajectories parts begun on the second type circles and ended in the fixed 
2-point. Then surface $F$ along with these disks cuts 3-manifold on 3-disks, 
each of which has one fixed point of index 0 or 3. Having the homeomorphisms of
the boundaries of these, shall continue their into interior of its. Thus we 
constructed homeomorphism of manifolds which set the topological equivalence of
vector fields.

\bigskip
{\bf 3. The expanding of the isomorphisms of graphs up to surfaces 
homeomorphism}. 
\medskip

Let $G$ is an oriented graph imbedded into surface $F$, and $G' $ into 
$F'$. If the graphs are isomorphic (it is possible not preserving the 
orientation of edges) then there exist the finite number of the different 
isomorphisms between them. Then question on the existence of surfaces 
homeomorphism, restriction of which on graphs is the graphs isomorphism, is 
equivalent to question on the capability of the graphs isomorphism extension up
to surfaces homeomorphism. 

Let $g: G \rightarrow  G'$ is a graphs isomorphism which maps vertex 
$A_i$ of the graph $G$ on vertex $A'_i$ , and edges $B_j$ on $B'_j$.  

Let $g: G \rightarrow  G' $ be a graph isomorphism which maps vertex 
$A_i$ of graph $G$ on vertex $A'_i$ , and edges $B_j$ on 
$B'_j$. Denote by $U(G)$ the tubular neighborhood of graph $G$ in 
surface $N$ and let  be the projection of its closure on \`a graph. Then 
complement $N\backslash U(G)$ consist of surfaces $F_i$ with boundary $and 
\cup \partial F_i= \partial U(G)$. Let us cut each circle from the boundary of 
surface $F_i$ to arcs in a such way that each arc maps by projection $p$ on one
edge of graph $G$ and reverse image $p-^1(B_j)$ for each edges consist of two 
arc from all surfaces. Let us choose such orientation of the arcs that the 
projection $p$ preserves the orientation and we will denote these arcs by the 
same letters as the appropriate edges. 

We fix the orientation on each surface $F_i$ (which is compatible with surface 
$F$ orientation if the surface $F$ is oriented and in an arbitrary way 
otherwise). For each circle from the boundary of surface we form a word 
consisted from letters $B^{\pm 1}_j$ which denote arcs (edges of graph) from 
this circle. We write the letters in such consequence in which we meet it when 
we go around the circle along orientation compatible with surface $F_i$ 
orientation. Letter has degree $+1$ if the orientation of correspondence arc is
the same as the circle orientation and -1 otherwise. Two words are called 
equivalent if one from another can be obtained by cyclic letters permutation. 
This situation is obtained if we choose the other beginning of the circuit. The
words are called reverse if one from another can be obtained in result of 
writing letters in reverse order with changing it degree and, possibly, cyclic 
permutation. This situation holds under the circle circuit with other 
orientation. 

For each surface $F_i$ we compose the list consisted of the number $n_i$ which 
is equal the surface $F_i$ genus and words which have written when we go around
surface boundary along orientation. Two such list are called equivalent if it 
has the same numbers $n_i$ and there is the one to one correspondence between 
words such that corespondent words are equivalent. The list are reverse if it 
has the same numbers $n_i$ and all the words are reverse. 

Thus for surface $N$ and of graph $G$ we construct the collection of lists in a
way that each list correspond one surface $F_i$. Two such collection are called
equivalent if there is one to one correspondence between list such that 
correspondent list are equivalent or reverse. 
\medskip
{\bf Lemma 2.} Let $G$ be the oriented graph embedded in surface $N$ and 
$G' $ in $N'$, $g: G \rightarrow  G'$ be a graph 
isomorphism, which maps vertex $A_i$ of graph $G$ in vertex $A' _i$ and
edges $B_j$ in $B'_j$. Then the graph isomorphism can be extended to 
surface homeomorphism iff replacing $B_j$ on $B'^{\pm 1}_j$ (we choose
sign in depending on edges orientation preserving) in the lists collection for 
pair $(N,G)$ we obtain the lists collection for pair 
$(N',G' )$. In addition the homeomorphism preserve the 
orientation if all correspondent lists are equivalent. 
\medskip
Proof. Necessity of theorem condition is followed from construction. Indeed 
surfaces homeomorphism gives one to one correspondence between the lists and 
the word sets in its, in addition if we get start go around from other point we
obtain equivalent words and lists, and if we reverse the orientation we obtain 
reverse ones.  

Sufficiency. Let $U_i$ be a connected component after cutting surface $N$ by 
graph. Then they are gomeomorphic to interiors of surfaces $F_i$. Surface $N$ 
can be obtained in result of surfaces $F_i$ gluing to graph G. Let us consider 
the boundaries of surfaces $F_i$ as $n$-tagon (where $n$ is the number of 
letters in the word). Each edge corresponds to one letter from word and 
attached map on one edge of graph G. Then the graph isomorphism and the word 
equivalence give natural homeomorphism between surface boundaries. Because of 
the genus and the number of boundary connected component are same for surfaces 
$F_i$ and $F_i' $ then they are gomeomorphic. In addition there exist 
homeomorphism, which expand the given boundary homeomorphism. This means that 
graph isomorphism can be expanded to surface homeomorphism.

\bigskip
{\bf 4. The gradient like vector field invariant construction.}

\medskip
As we do it in section 2 for each vector field we construct surface with two 
sets of circles on it. The graph is the union of these circles. Vertexes of 
graph are the circle intersection and one arbitrary point on each circle 
without intersections. The edges are the arc between vertexes. All edges of 
graph are decomposed on two sets depending on which sets of circle the 
correspondent circle belong. We fix the arbitrary orientation of this graph. 
For this embedded graph as in section 3 we construct the word list collection 
with letters corresponded to the edges of the graph.  

Definition. Such constructed graph with edges decomposition on two sets and 
with word list collection is called distinguished graph of vector field. Two 
distinguished graphs are called equivalent if there exist graph isomorphism, 
which preserve edges decomposition on two set. Replacing letters from word list
collection of first graph by corresponding letters from second graph (with 
degree $\pm 1$ depending on orientation) we have to obtain word list 
collection, which is equivalent to second graph word list collection. 

{\bf Theorem 1.} Two vector fields are topological equivalent if an only if 
their distinguished graphs are equivalent.  

Proof is followed from using lemmas 1 and 2.

\medskip
{\bf 5. Topological conjugation of diffeomorphisms.} 
\medskip
Let $f: M^3\rightarrow  M^3$ be a gradientlike Morse-Smale diffeomorphism. As 
it have been done in section 2 we construct surface with two set of circle on 
it and afterwards as in section 4 distinguished graph. Then the diffeomorphism 
$f$ action on saddle point integral manifolds induce map between first type 
circles and map between second type circles and isomorphism of distinguished 
graph on itself, which we call by inner.  

{\bf Theorem 2.} Two gradientlike Morse-Smale diffeomorphisms $f$ and $g$ are 
topological conjugated iff there exist their distinguished graphs isomorphism 
which gives the equivalence between them and commutate with inner graph 
isomorphisms. 

Proof. Necessity follow from construction. Let us prove the sufficiency. As in 
theorem 1 we construct homeomorphism $h$ between manifolds which maps stable 
manifold on stable ones and unstable on unstable. In addition $h(f (U)) = g 
(h(U))$, where $U$ is the part of stable or unstable manifold on which it 
separated by another manifolds. Analogously dimension 2 in [Grines] this 
homeomorphism can be corrected to needed homeomorphism.

\bigskip
{\bf 6. Realization of dynamic system with given invariant.}

\medskip
Let us research the problem when distinguishing graph represent a surface with 
two sets of circles on them and gradient like vector field. Let $K$ is a 
complex received in result of gluing surfaces, appropriate to the lists, to the
graph. As each list of words sets a surface with boundary, structure of a 
surface (to be local homeomorphic to the plane) can be broken only in gluing 
points, that is on edges and in vertexes of the praph. 

\noindent
1) The condition that, that a complex $K$ is locally plane in internal points 
of edges equivalent to that there is two part of surfaces boundaries that is 
glued to the edge. It means, that each letter or it reverses meets equally two 
times in all lists. 

\noindent
2) We shall assume that the first condition is executed. For each vertex we 
shall consider the set of incident edges. Two edges we shall name adjacent if 
they lie in the boundary of one of gluing surfaces and have the common vertex 
in it. This condition equivalently to that there is a word in which the 
appropriate letters are adjacent or are first and last letters of a word. Two 
edges we shall name equivalent if there is a chain of the adjacent edges, which
connect them. Then a condition that that a complex $K$ is local plane in vertex
is equivalent to that for each vertex all incident edges are equivalent. 

\noindent
3) Let conditions 1) and 2) are executed. Then distinguishing graph determinate
the graph on a surface. This graph is two sets of circles if and only if each 
vertex is incident to four edges or is the beginning and end of one edge, if it
together with this edge forms a loop. Further, in vertex with four incident 
edges the adjacent edges should belong to various sets of edges (which 
correspond to two sets of circles). 

\medskip
{\bf Theorem 3}. Distinguishing graph is a graph of gradient-like vector field 
if and only if

a) Each letter in the whole set of the lists of words meets equally two 
times, 

b) Each vertex is incident to four edges or is the beginning and end of one 
edge,

c) Each of four edges, is adjacent for two others from other set of edges.

\medskip
Proof is followed from discussion above.  

\bigskip
{\bf Literature.}

\medskip
1. Aranson S.H., Grines V.Z. Topological classification of flows on closed 
2-manifolds $//$ Advant. math sc.(Russian)- 41, N.1, 1986.- P.149-169. 

2. A.V.Bolsinov, A.A.Oshemkov, V.V.Sharko. On Classification of Flows on 
Manifolds. I $//$ Methods of Funct. An. And Topology, $v.2$, N.2, 1996, 
$P.131-146.$ 

3. Peixoto M. On the classification of flows on two-manifolds. In: Dynamical 
Systems, edited by M.Peixoto.- Academic Press.- 1973.- P.389-419.  

4. Fleitas G. Classification of Gradient like flows of dimensions two and three
$//$ Bol. Soc. Brasil. Mat.-9, 1975.- P.155-183. 

5. Umanski Ya.L. The circuit of three-dimensional Morse-Smale dynamic system 
without the closed trajectories $//$Isvestiya of USSR Acad. of Sc., 230, N6,1976, $P.1286-1289.$ 

6. Grines V.Z. Kalay I.N. Classification of three-dimensional gradient like 
dynamic systems $//$ Advant. Math. Sc.(Russian), $v.49$, N.2, 1994.- P.149-150. 

7. Smale. S. On gradient dynamical system $//$ Ann. Math. 74, 1961, P.199-206. 

8. Anosov D.V. Smooth dynamic systems 1 $//$ Results of a science and 
engineering. Modern math. Problems. Fund. Directions, V.1, 1985,(Russian)- P.151-242.

9. Matveeev S.V., Fomenko A.T. Algorithmic and computer methods in 
three-dimensional topology. --- M.: MSU, 1991.- 301p. 

10.Prishlyak A. On graph, which is embedded in a surface $//$ Advant. Math. Sc. 
(Russian), V.56, N.4, 1997, -P.211-212.

\medskip
Kiev University

\medskip
e-mail: prish@mechmat.univ.kiev.ua
\vfill\eject\end